\documentclass[conference]{IEEEtran}
\ifCLASSINFOpdf
  \usepackage[pdftex]{graphicx}
  \graphicspath{{./figures/}}
\else
\fi
\usepackage{amsmath}
\usepackage{algorithmic}
  \usepackage[caption=false,font=footnotesize]{subfig}
\usepackage{url}

\usepackage{amsmath}
\usepackage{amsfonts}
\usepackage{amssymb}

\usepackage{accents}

\newcommand{\R}{\mathbb{R}}
\newcommand{\B}{\mathbb{B}}
\newcommand{\N}{\mathbb{N}}

\newcommand{\tl}[1]{\tilde{#1}}
\newcommand{\ubar}[1]{\underaccent{\bar}{#1}}

\def\cL{\mathcal{L}}

\def\ij{{ij}}
\def\cC{\mathcal{C}_\ij}
\def\cT{\mathcal{T}}

\def\cLe{\cL_{\text{ext}}}

\begin{document}
%
\title{Heuristics for Transmission Expansion Planning 
in Low-Carbon Energy System Models}

\author{\IEEEauthorblockN{Fabian Neumann and Tom Brown}
\IEEEauthorblockA{Institute for Automation and Applied Informatics\\
Karlsruhe Institute of Technology\\
fabian.neumann@kit.edu
}
}


%


\IEEEoverridecommandlockouts
\IEEEpubid{\makebox[\columnwidth]{978-1-7281-1257-2/19/\$31.00~\copyright2019 IEEE \hfill} \hspace{\columnsep}\makebox[\columnwidth]{ }}

\maketitle

\IEEEpubidadjcol

\begin{abstract}

  Governments across the world are planning to increase the share of renewables in their energy systems.
  The siting of new wind and solar power plants requires close coordination with grid planning,
  and hence co-optimization of investment in generation and transmission expansion
  in spatially and temporally resolved models is an indispensable but complex problem.
 Particularly considerations of transmission expansion planning (TEP) add to the problem's complexity. Even if the power flow equations are linearized, the optimization problem is still bilinear and mixed-integer due to the dependence of line expansion on line impedance and a discrete set of line expansion options. While it is possible to linearize this mixed-integer nonlinear program (MINLP) by applying a big-$M$ disjunctive relaxation, the resulting MILP is still hard to solve using state-of-the-art solvers for large-scale energy system models.

  In this paper we therefore develop heuristics to incorporate integer transmission expansion and
  responsive line impedances in energy system models,
  while retaining the lower computational effort of continuous linear programming (LP)
  by applying sequential linear programming (SLP) techniques, relaxation and discretization approaches.
  We benchmark their performance against the results of the exact formulation
  for a policy-relevant case study of the German transmission system
  in terms of their speed-up in computation time,
  deviation from optimal total system cost, and similarity of line expansion. 
  
  Using heuristics we reduce computation times of joint generation and transmission optimisation by 82\% with a maximal total system cost deviation of only 1.5\%.
  The heuristics closely mirror optimal integer line investment of the exact MINLP with 
  considerable time savings for policy-relevant low-carbon energy system optimization models.

\end{abstract}


%
\IEEEpeerreviewmaketitle

\section{Introduction}

Governments across the world are planning to increase the share of renewables in their energy systems.
The siting of new wind and solar power plants requires close coordination with grid planning.
Co-optimization of generation and transmission expansion is indispensable
to weigh the use of local renewable resources against the costs of grid extension to sites with better resources \cite{Krishnan2016}.
While co-optimization is not always possible within
the regulatory structures of some markets, the impetus may grow as
public acceptance limits how much grid expansion is possible. 
However, co-optimization demands high spatial detail to factor in grid bottlenecks and renewable potentials,
and temporal detail to capture the variability of renewables and loads \cite{Hoersch2017d,Kotzur2018a}.
Together, this leads to a high problem complexity. 
Therefore, as a necessary reduction in accuracy, many energy system models either employ a simple transport model, which ignores grid physics,
or they assume a linearized power flow formulation but don't account correctly for how grid impedances change as lines are expanded \cite{Groissbock2019, Lumbreras2016}.

On the other hand, studies focusing more specifically on large-scale transmission expansion planning (TEP), i.e. without generation planning, 
typically model the available line types and the change of power flows in response to line expansion more realistically \cite{Lumbreras2017,Lumbreras2017a}.
While these studies are capable of handling multi-thousand node networks if the considered range of operating scenarios is severely limited,
they usually assume fixed generator capacities and thereby disregard the investment trade-offs
revealed by the joint optimization of transmission and generation capacities.

Merging the two domains entails a complex optimization problem.
Besides many integer variables for line expansion, the problem includes highly non-linear and non-convex power flow.
Even if the power flow equations are linearized, the optimization problem is still bilinear and mixed-integer
due to the dependence of line expansion on line impedance \cite{Krishnan2016, Hagspiel2014}.
One can further linearize this problem by applying a \mbox{big-$M$} reformulation to obtain a mixed-integer linear program (MILP) \cite{Krishnan2016}.
And yet, with an increasing size of power networks, temporal resolution, and number of line expansion options,
even solving this reformulation with state-of-the-art algorithms becomes immensely time-consuming.

In this paper we therefore develop heuristics to incorporate integer transmission expansion and responsive line impedances into energy system models, while retaining the lower computational effort of continuous linear optimization problems (LP) by applying sequential linear programming (SLP) techniques.
An LP, that jointly optimises generation and transmission capacities, is solved iteratively.
In each iteration, the line impedances are adjusted to the current line capacities.
Subsequently, the line extensions are fixed to their closest integral values followed by a final iteration of generation expansion.
Such heuristics are of large avail when e.g. in a model with storage expansion the merits of frequently applied Benders or other decompositions cannot be leveraged because individual snapshots can no longer be decoupled and the resulting Benders cuts become ineffective \cite{Ramos2016}.

SLP for joint optimization of transmission and generation
was introduced in \cite{Hagspiel2014} but it was only tested against the exact
solution for a 3-node showcase example and the authors did not benchmark different
variants of the discretisation against the exact solution for realistic
case studies.

Here we remedy these gaps in the
literature by evaluating the performance of different heuristic versions of
SLP for a policy-relevant case study of the German transmission system.
We compare the heuristics in terms of their speed-up in computation time, deviation from optimal total system cost, and similarity of line expansion. 
Furthermore, we show the error made by disregarding integrality constraints and responsive line impedances.

The remainder of the paper is organized as follows: 
Section \ref{sec:problem} guides through the mathematical formulation of the generation and transmission expansion planning problem and describes its standard solving algorithms, while Section \ref{sec:heuristics} elaborates on heuristic algorithms for the problem at hand.
Section \ref{sec:case} describes our case-study, the results of which are discussed in Section \ref{sec:results}.
The work is critically appraised in Section \ref{sec:appraisal} and concluded in Section \ref{sec:conclusion}.

\section{Problem Formulation}
\label{sec:problem}

The techno-economic objective is to minimize the total annual system costs, composed of annualized capital costs $c_*$ for generator capacity $G_{i,s}$ of technology $s$ as well as line capacities $F_{ij}$, HVDC link capacities $H_{ij}$ and variable operating costs $o_*$ for generator dispatch $g_{i,s,t}$:

\begin{equation}
  \min \quad \sum_{ij} c_{ij}^L F_{ij} + \sum_{ij} c_{ij}^H H_{ij} + \sum_{i,s} c_{i,s} G_{i,s} + \sum_{i,s,t} w_t o_{i,s} g_{i,s,t}
\end{equation}
where representative points in time are weighted by $w_t$ such that $\sum_t w_t = 365\cdot 24 \text{h} = 8760 \text{h}$.

The objective function is subject to a set of linear constraints $\mathcal{S}$ specifying limits on 
(i) the capacities of generators, lines and links from geographical and technical potentials, 
(ii) the availability of renewables per snapshot and location derived from re-analysis weather data, 
(iii) the minimum share of renewable energy generation, and
(iv) the nodal balance equations as demanded by Kirchhoff's current law.
Equations for Kirchhoff's voltage law (v), complete the linearized optimal power flow (LOPF) formulation, also known as DC-approximation  \cite{Brown2018}.

The flow $f_{ij,t}$ of the transmission line connecting buses $i$ and $j$ at snapshot $t\in\cT$ is

\begin{equation}
  f_{ij,t} = b_{ij} (\theta_{i,t} - \theta_{j,t})
  \label{eq:flows}
\end{equation}
where $b_{ij}$ is the susceptance and  $\theta_{i,t}$ denotes the voltage angle.

It is a common assumption in current energy system models that line susceptances are invariant to line investments ($b_{ij} = \tl b_{ij}$)
and in some models line capacity is expanded continuously.
On the other hand, in addition to acknowledging that the susceptance changes proportional to the capacity of a line ($x_{ij}^{-1}=b_{ij}\propto F_{ij}$),
TEP studies address the limited choice of circuit types.
Instead of allowing the transmission capacities to expand continuously, they consider the number of added circuits $\Gamma_{ij}\in \N$ to be within a set of candidate circuits $\cC \subset \N_{\geq 0}$ for every extendable line $ij\in \cLe$.

The line capacity $F_{ij}$ is linked to the choice of $\Gamma_{ij}$ via

\begin{equation}
  F_{ij} = \left(1+\frac{\Gamma_{ij}}{\tilde \gamma_{ij}}\right)  \tilde F_{ij} \label{eq:capacity}
\end{equation}
where $\tl \gamma_{ij}\in\N$ denotes the initial number of circuits and $\tl F_{ij}$ is the initial line capacity. 

Amending equation (\ref{eq:flows}) with adaptive impedances and lumpy investment entails a set of non-linear and non-convex constraints

\begin{align}
  f_{ij,t} = \overbrace{\left(1+\frac{\Gamma_{ij}}{\tilde \gamma_{ij}}\right)  \tilde b_{ij}}^{b_{ij}} \cdot (\theta_{i,t} - \theta_{j,t}) \label{eq:adflows}\\
  \Gamma_{ij}\in \cC \subset \N_{\geq 0} \label{eq:intinv}
\end{align}
where $\tl b_{ij}$ is the initial line susceptance.

With equations (\ref{eq:adflows})-(\ref{eq:intinv}) the optimization problem classifies as a hard-to-solve mixed-integer nonlinear problem (MINLP) for which this paper seeks to find heuristic solutions. 

\subsection{Big-$M$ Disjunctive Relaxation}

However, equations (\ref{eq:adflows})-(\ref{eq:intinv}) can be reformulated into a mixed-integer linear problem (MILP) using a \mbox{big-$M$} disjunctive relaxation at the cost of numerous extra variables and constraints \cite{Krishnan2016}.
The \mbox{big-$M$} reformulation uses one or multiple large constants $M$ in combination with binary variables to formulate disjunctive linear inequalities that redefine originally nonlinear constraints.

For the given problem, the integer variables $\Gamma_{ij}\in \cC \subset \N_{\geq 0}$ are substituted with binary variables $\Gamma_{ij,c}\in\B$ for each candidate investment $c\in\cC$ and the constraint

\begin{equation}
  \sum_{c\in\cC} \Gamma_{ij,c} = 1
\end{equation}
is added to allow only one candidate choice per line. Note that for all lines no capacity expansion is a valid choice, i.e. $\forall ij \in \cLe:0\in\cC$.
The resulting disjunctive inequalities for all line investment candidates $c\in\cC$ 

\begin{align}
  \left(1+\frac{c}{\tl \gamma_\ij}\right) \tl b_{ij} \left( \theta_{i,t} -  \theta_{j,t}\right) -  f_{\ij,t} \geq ( \Gamma_{\ij,c} - 1 ) \bar{M}_\ij \label{eq:kvl1} \\
  \left(1+\frac{c}{\tl \gamma_\ij}\right) \tl b_{ij} \left( \theta_{i,t} -  \theta_{j,t}\right) -  f_{\ij,t} \leq (1 -  \Gamma_{\ij,c}) \ubar{M}_\ij \label{eq:kvl2}
\end{align}
are in effect equivalent to the equality constraint (\ref{eq:adflows}) if the \mbox{big-$M$} parameters $M_\ij$ are suitably chosen.
If $\Gamma_{ij,c}$ is $0$, the inequalities (\ref{eq:kvl1})-(\ref{eq:kvl2}) are inactive, but if $\Gamma_{ij,c}$ is $1$, the right-hand side is $0$ and both merge to an equality constraint.

The choice of $M_\ij$, however, is non-trivial.
While the values must be chosen as large as necessary to guarantee
that for a certain investment decision all other constraints linked to a different investment decision are inactive,
they should also be as small as possible to avoid numerical instabilities.
One approach for determining minimal values for $M_\ij$ using shortest-path optimization is outlined in \cite{Binato2001}.

Adjusting equation (\ref{eq:capacity}) to the binary nature of the line investment variables completes the reformulation.

\begin{equation}
  F_{ij} = \left(1+\frac{
    \sum_{c\in\cC} c \cdot \Gamma_{ij,c}
  }{
    \tilde \gamma_{ij}
  }\right)  \tilde F_{ij}
\end{equation}

While applying the \mbox{big-$M$} formulation resolves the nonlinearities of equation (\ref{eq:adflows}),
the reformulation retains the integrality of line investment and might therefore still be hard solve quickly.
More precisely, the reformulation considerably adds to the problem size in terms of the number of additional variables and constraints, where the cardinality of a set is denoted by $|\cdot|$:

\begin{align}
  \text{variables:}\quad &\sum_{ij\in\cLe} |\cC| - |\cLe|\\
  \text{constraints:}\quad &|\cLe| + |\cT| \cdot \left(\sum_{ij\in\cLe} |\cC| - |\cLe|\right)
\end{align}

In this regard, the combinatorial difficulty of TEP has led researchers to focus only on promising candidate investments; e.g. an automatic candidate selection scheme is discussed in \cite{Lumbreras2014}.

\subsection{Benders Decomposition}

The computational challenges have led many researchers to apply decomposition techniques such as Benders decomposition (BD) to TEP (e.g. \cite{Lumbreras2017}). It is therefore worth briefly discussing its suitability for the problem at hand.
BD divides the optimization problem into an investment master problem with integer variables and
an optimal power flow subproblem with exclusively continuous variables.
By iteratively adding cuts to the master problem, derived from the duals of the subproblem given a specific set of investments,
the operational expenses are approximated \cite{Binato2001}.

Despite its prevelance in the literature, researchers have raised concerns about the scalability of a straightforwardly applied BD algorithm for large-scale TEP problems
leading to investigations on acceleration techniques \cite{Majidi-Qadikolai2018, Munoz2017}.
In \cite{Ramos2016}, a speedup by around factor 2 was attained by
generating cuts based on relaxations and suboptimal solutions,
removing inactive cuts,
and, most significantly, by applying multicuts.
The subproblem is split into individual snapshots and each of the resulting subproblems adds its own cut to the master problem.
However, the merits of multicuts can only be leveraged if snapshots are independent of each other.

However, in the presence of storage and renewable energy targets, there is a significant degree of intertemporal coupling.
Additionally, full energy system models involve significantly more decision variables
and, thereby, many more trade-offs than pure TEP studies.
Hence, while the value of BD for detailed TEP studies remains unquestioned,
for energy system models working with heuristics and avoiding the NP-hardness of integer problems could prove more suitable.

\subsection{Tuned Parameter Set for Solver Settings}

In this study we therefore solve the big-$M$ formulation solely with a tuned parameter set listed in Table \ref{tab:gurobi} for the solver Gurobi. The MIP optimality gap (MIPGap) is a parameter which relates the lower and upper bounds of the objective value throughout the optimization process and serves as termination condition.

\begin{table}
  \caption{Solver settings for solving MILP in Gurobi}
  \begin{center}
    
    \begin{tabular}{ll}
      \hline
      \textbf{Parameter} & \textbf{Value} \\
      \hline
      MIPGap & 0.015, 0.01, 0.005 \\
      MIPFocus & 2 \\
      Heuristics & 0.5 \\
      Cuts & 1 \\
      NumericFocus & 3 \\
      TimeLimit & 72 hours \\
      \hline
    \end{tabular}
  \end{center}
  \label{tab:gurobi}
\end{table}

\section{Heuristic Algorithms}
\label{sec:heuristics}

The proposed heuristics in this paper comprise three central elements:
relaxation of integer line investment (Section \ref{sec:relaxation}),
iterative updates of line impedances (Section \ref{sec:iteration}), and
post-facto discretization of line investment (Section \ref{sec:discretisation}).
We label the different heuristic approaches presented above with abbreviations listed in Table \ref{tab:experiments}.

\begin{table}
  \caption{Experiment code glossary}
  \begin{center}
    
    \begin{tabular}{ll}
      \hline
      \textbf{Code} & \textbf{Explanation} \\
      \hline
      heur & Heuristic approach to solve MINLP \\
      int & Integer line investment \\
      iter & Iterating until convergence  \\
      seqdisc & Sequential discretization of impedances \\
      postdisc & Post-facto discretization with single threshold \\
      postdisc-mult & Post-facto discretization with multiple thresholds \\
      \hline
    \end{tabular}
  \end{center}
  \label{tab:experiments}
\end{table}

\subsection{Relaxation of Line Investment Variables}
\label{sec:relaxation}

To avail of the computational merits of continuous linear programs,
the integer investment decisions $\Gamma_{ij} \in \cC\subset \N_{\geq 0}$ (tagged \textit{heur-int})
may be relaxed to allow each line to be expanded continuously, i.e. $\Gamma_{ij} \in \mathbb{R}$ (tagged \textit{heur}).

\subsection{Iterative Update of Line Impedances}
\label{sec:iteration}


To resolve the problem's nonlinearities, we draw on the concept of sequential or successive linear programming (SLP).
In SLP, constraints are iteratively linearized by a first-order Taylor series expansion around an operating point
(e.g. currently optimal line capacity) \cite{Griffith1961}.
This way, instead of one complex nonlinear problem, multiple simpler linear programs are solved step by step.
It can be shown that when SLP converges, the Karush-Kuhn-Tucker (KKT) conditions\footnote{The KKT conditions are the first order necessary conditions for a solution of a nonlinear optimization problem to be optimal.} are satisfied \cite{Bazaraa2006}. 
Additionally, move limits, which constrain the variable variation between two subsequent iterations to a linear trust region,
are frequently supplemented but can be omitted if convergence is obtained without them \cite{Bazaraa2006}.

For the given problem, equation (\ref{eq:adflows}) is modified to

\begin{equation}
  f_{ij,t} = b_{ij}^{(k)} (\theta_{i,t} - \theta_{j,t})
\end{equation}
where in the first iteration the initial susceptances are used

\begin{equation}
  b_{ij}^{(1)} = \tilde b_{ij}
\end{equation}
while for subsequent iterations $k+1$ the optimal line investment $\Gamma_{ij}^{*(k)}$ of the previous iteration $k$ determines the physical line characteristics

\begin{equation}
  b_{ij}^{(k+1)} = \left(1+\frac{\Gamma_{ij}^{*(k)}}{\tilde \gamma_{ij}}\right)\tilde b_{ij} \qquad \forall k>1
\end{equation}

Instead of adjusting the susceptances to values corresponding to fractional line capacities (tagged \textit{iter}),
another variant is to round any $\Gamma_{ij}^{*(k)}\in\R$ to their nearest integer value,
referred to as sequential impedance discretization (tagged \textit{iter-seqdisc}).

The iteration loop terminates either if a pre-defined iteration limit is reached or if there is convergence; i.e. 
in two consecutive iterations there is no change in line investment or objective value and, therefore, also $b_{ij}^{(k+1)} = b_{ij}^{(k)}$.
While in \cite{Hagspiel2014} only a loose convergence tolerance for the objective function of 500 MEUR/a was required, we tighten it to 1000 EUR/a.

\subsection{Post-facto Discretization of Line Investment Variables}
\label{sec:discretisation}

Since optimal line capacities are likely to be fractional when the continuous relaxation is applied,
they do not represent a valid investment choice according to the candidate set $\cC$.
Post-facto discretization follows the iteration loop to produce a valid set of investment decisions.

The optimal capacities of continuously expanded lines are directly rounded to an integer value
using a discretization threshold ($z=0.3$) on the fractional investment (tagged \textit{postdisc}) \cite{Chekuri2009}.
The discretized line capacities are then fixed and line impedances are adapted accordingly for a final round of generation expansion only, 
to find the optimal mix of generator and link capacities given the discretized line capacities.

A more extensive post-facto discretization procedure (tagged \textit{postdisc\_mult}) repeats the steps above for multiple discretization thresholds ($z \in \{0.1,0.2,0.3,0.4,0.5\}$)
and chooses the configuration entailing the lowest total system costs.
While this approach bears the chance of outcomes closer to the global optimum, it will naturally take longer to solve if thresholds are serially evaluated.
However, this step lends itself to parallelization to reduce computation times.

The post-facto discretization step ensures that, ultimately, a feasible solution of the original problem is obtained.\\

\section{Case Study: TEP in the German Power System}
\label{sec:case}

The presented transmission expansion heuristics are evaluated on simplified models of the German power sector,
because they exhibit a distinct spatial mismatch of abundant wind resources in the North 
and high loads in the South, suggesting high future strains on the transmission network.
Already today (2017), grid congestion causes curtailment of 5.52 TWh corresponding to 2.9 \% of renewable energy produced, and has been growing rapidly \cite{Bundesnetzagentur2018}. 

\subsection{Model Building and Assumptions}
\label{sec:model}

The case-study model is obtained using PyPSA-Eur, an open model dataset of the European power system comprising the ENTSO-E transmission network.
Full details on the underlying dataset and assumptions can be found in \cite{Horsch2018b};
complementary considerations are outlined below.
We allow simultaneous capacity expansion of transmission lines, HVDC links, and various types of generators:
photovoltaic, onshore wind, and offshore wind generators with AC and DC grid connections, as well as open- and combined-cycle gas turbines.
Currently existing run-of-river and biomass capacities are not extendable.

All transmission lines are assumed to be of the same type `Al/St 240/40 4-bundle 380.0' with an Aluminium/Steel cross-section of 240/40 and a 4-bundle of wires per phase at \mbox{380 kV} and 
allow for two additional parallel circuits fitted to the existing towers \cite{Oeding2011}. 
To avoid issues with social acceptance, a limit on the total additional volume of transmission capacities of 25\% measured in MWkm is imposed.
Moreover, to approximate $N-1$ security and withhold capacity for reactive power flows,
line capacities may not be loaded by more than 70\% of their nominal power rating \cite{Horsch2018b}.

The corridors for HVDC links are taken from the TYNDP 2018 and may be expanded continuously with a capacity restriction of 8 GW each.
Stub links for offshore wind power plants with DC grid connection are exempt from this capacity restriction.

We neglect storage options (e.g. pumped-hydro, batteries, and power-to-hydrogen) to further simplify the case study and direct the focus to transmission expansion.
Since balancing renewables in time through storage is an imperfect substitute for balancing renewables in space through transmission networks,
we note that ignoring storage options rather overestimates the optimal grid reinforcements.
A target share of renewable electricity of 70\% is chosen since it is high enough to require transmission expansion and low enough to be viable without storage options.
Current policy aims at a 65\% share of renewable electricity by 2030.

To assess the scaling of the heuristic algorithms both in terms of spatial and temporal resolution,
models with different levels of regional and time series aggregation are analyzed.
The German power network is reduced to $n\in\{20,40,60,80,100\}$ nodes using the network clustering functionality of the power system analysis toolbox PyPSA \cite{Brown2018},
while the number of hourly snapshots is reduced to $t\in\{100,200,300,400\}$ using the time series aggregation module \textit{tsam} \cite{Kotzur2018a}.
The cluster weights represent the snapshot weightings and add up to 8760 hours to match the annuities of capital expenditures.
Both reduced network size and temporal resolution are necessary limitations in accuracy to be able to solve the exact TEP formulations of Section \ref{sec:problem} for up to 200 snapshots against which the heuristics of Section \ref{sec:heuristics} are evaluated and may not accurately represent the necessary range of operating conditions to make qualified investment decisions. 

\subsection{Experimental Setup}

The experimental setup comprises solving the case-study with the ensemble of network sizes and snapshots outlined in Section \ref{sec:model} for the exact benchmark formulation that uses the big-$M$ formulation described in Section \ref{sec:problem}
and the different variants of the proposed heuristic transmission expansion planning approaches named by combinations of the codes listed in Table \ref{tab:experiments} and described in Section \ref{sec:heuristics}.
Note, that the approaches \textit{heur} and \textit{heur-iter} do not involve any discretization and the former neither considers changing impedances. Therefore, they do not provide feasible solutions, but represent current practice in some energy models.
The solution to the big-$M$ formulation comprises the lower and upper bounds that could be obtained for the lowest of the MIP optimality gap tolerances of Table \ref{tab:gurobi} within a walltime of 72 hours. It is further worth noting that the upper bound of the big-$M$ formulation is the first feasible solution the solver could find that satisfies this tolerance.


\section{Results and Discussion}
\label{sec:results}

\begin{figure*}
  \centering
  \includegraphics[width=\textwidth, trim=0cm 0.4cm 0cm 0.2cm, clip]{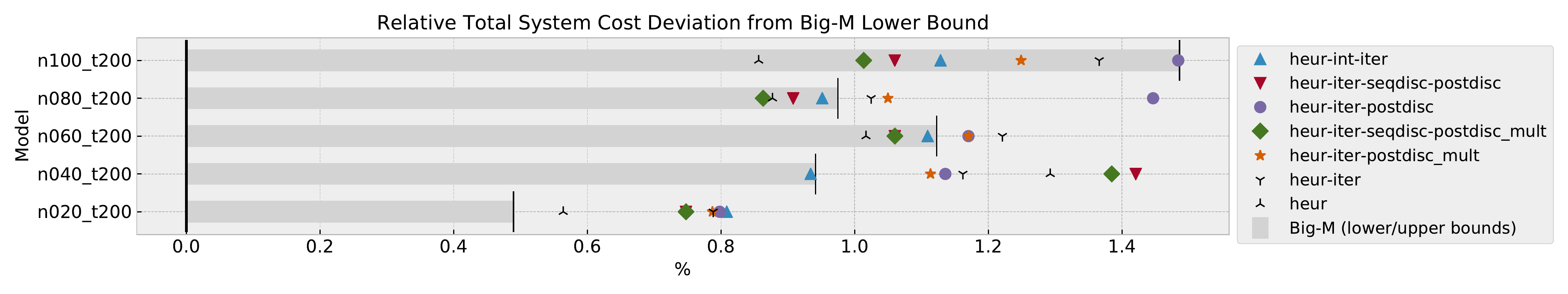}
  \includegraphics[width=\textwidth, trim=0cm 0.4cm 0cm 0.2cm, clip]{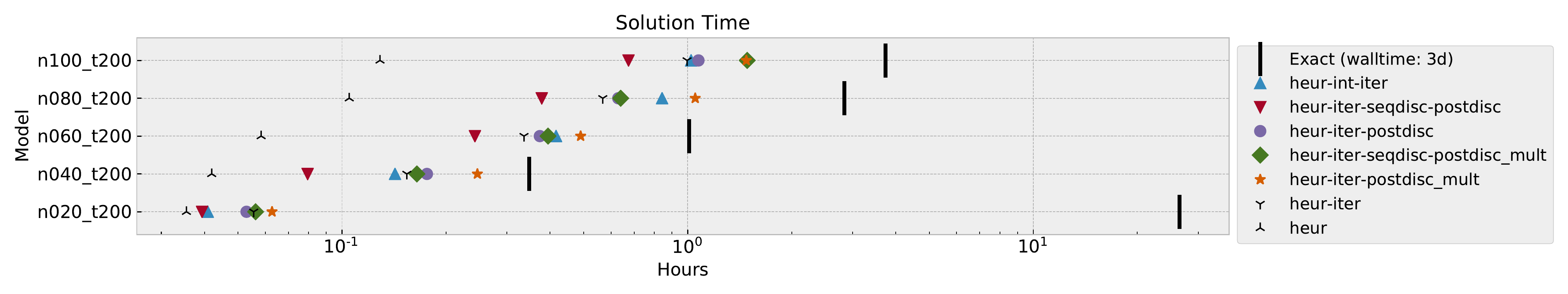}
  \includegraphics[width=\textwidth, trim=0cm 0.4cm 0cm 0.2cm, clip]{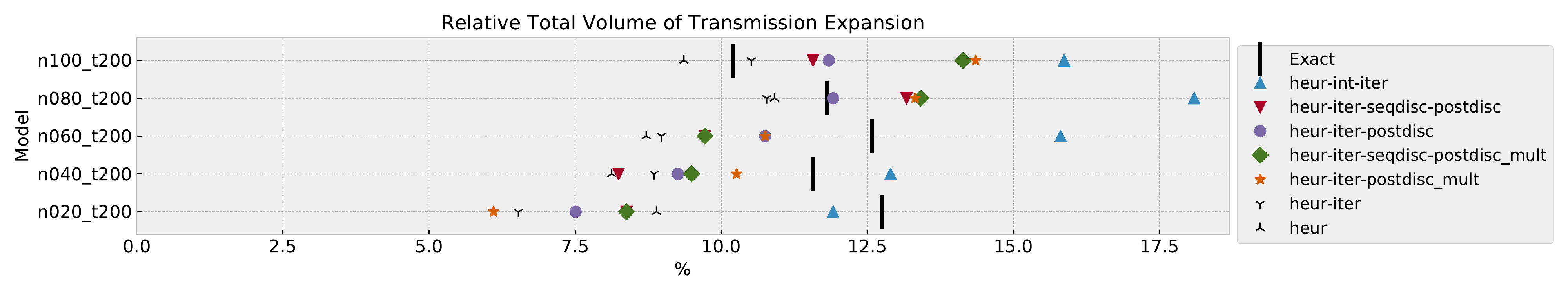}
  \vspace{-0.25cm} \caption{Comparison of 
  (i) relative total system cost deviation from the lower objective bound of the big-$M$ formulation,
  (ii) solving times, and 
  (iii) relative volume of added transmission lines over the volume of original transmission lines (TWkm) 
  across all presented heuristic approaches (as labelled in Table \ref{tab:experiments}) and the upper objective bound of the big-$M$ formulation (a feasible solution)}
  \label{fig:costtimevolume}
\end{figure*}

We assess the heuristics in terms of
(i) their deviation from optimal total annual system costs,
(ii) the solving time, and 
(iii) the similarity of line expansion in relation to the solutions of the exact big-$M$ reformulation.
Figure \ref{fig:costtimevolume} depicts these evaluation criteria for the network sizes outlined in Section \ref{sec:model} with 200 snapshots.
The colored markers represent feasible solutions to the exact MINLP, whereas black markers depict solutions of continuous relaxations of line investment and are generally not feasible MINLP solutions.
The shaded area between the two vertical lines represents the range within which the optimal solution of the big-$M$ formulation is contained.

In terms of total system costs, all heuristics across all tested models are contained within a 1.5\% relative cost increase compared to the lower bound of the big-$M$ formulation obtained with the specified MIP optimality gap tolerances.
Since the optimal solution is located between the lower and upper bounds, the heuristics are therefore at most 1.5\% more expensive than the optimal solution.
Depending on the model, the spread among the heuristics lies between 0.1\% and 0.5\%, while a consistent observation is that the heuristics are rather close to the upper bound of the big-$M$ formulation.
In abolute numbers, total system costs roughly evolve around 32 billion EUR/a, of which around 2.5 billion EUR/a are attributed to transmission infrastructure. Thus, a relative cost deviation of around 1\% corresponds to an abolute deviation of about 320 million EUR/a.

In terms of solving time, not iterating to update line impedances \textit{(heur)} is naturally a fast, but infeasible approach.
On the other hand, finding a solution to the big-$M$ formulation consumes the most time.
On average, the slowest heuristic is already faster by a factor of 2.2, but the approach \textit{heur-iter-seqdisc-postdisc} is invariably the fastest among the MINLP-feasible heuristics with an average speed-up by factor 5.4.
This is due to the sequential impedance discretization \textit{(iter-seqdisc)} which causes the iteration process to converge on average already after 4.4 iterations,
 while heuristics without sequential impedance discretisation consistently required the maximum number of 10 iterations (cf. Table \ref{tab:iterations}).

\begin{table}
  \caption{Average number of iterations for heuristic approaches}
  \begin{center}
    
    \begin{tabular}{lr}
      \hline
     \textbf{Algorithm} & \textbf{Average \# Iterations} \\
     \hline
     heur-int-iter      &                5.1 \\
      heur-iter-postdisc   &             10.0 \\
      heur-iter-postdisc-mult &           10.0 \\
      heur-iter-seqdisc-postdisc(-mult) &        4.4 \\
      \hline
    \end{tabular}
  \end{center}
  \label{tab:iterations}
\end{table}

A wide spread among the approaches can be observed in terms of the total volume of transmission expansion,
which we measure with a product of line capacity (MW) and length (km) and relate it to the original transmission network.
The differences range between 5\% and 7\%. While \textit{heur-int-iter} leads to the most line expansion for all models,
it does not incur the highest total system costs.
Approaches without discretization (\textit{heur} and \textit{heur-iter}) can tailor line expansion more accurately to needs and,
thus, entail less transmission expansion, but do not produce MINLP-feasible solutions.


\begin{figure*}
  \begin{center}
    \begin{scriptsize}
      \setlength\tabcolsep{0pt}
      \begin{tabular}{ccc}
        (i) Big-$M$ formulation & (ii) heur-iter-seqdisc-postdisc & (iii) heur-iter \\
        \includegraphics[width=0.34\textwidth, trim=0cm 3.7cm 0cm 1cm, clip]{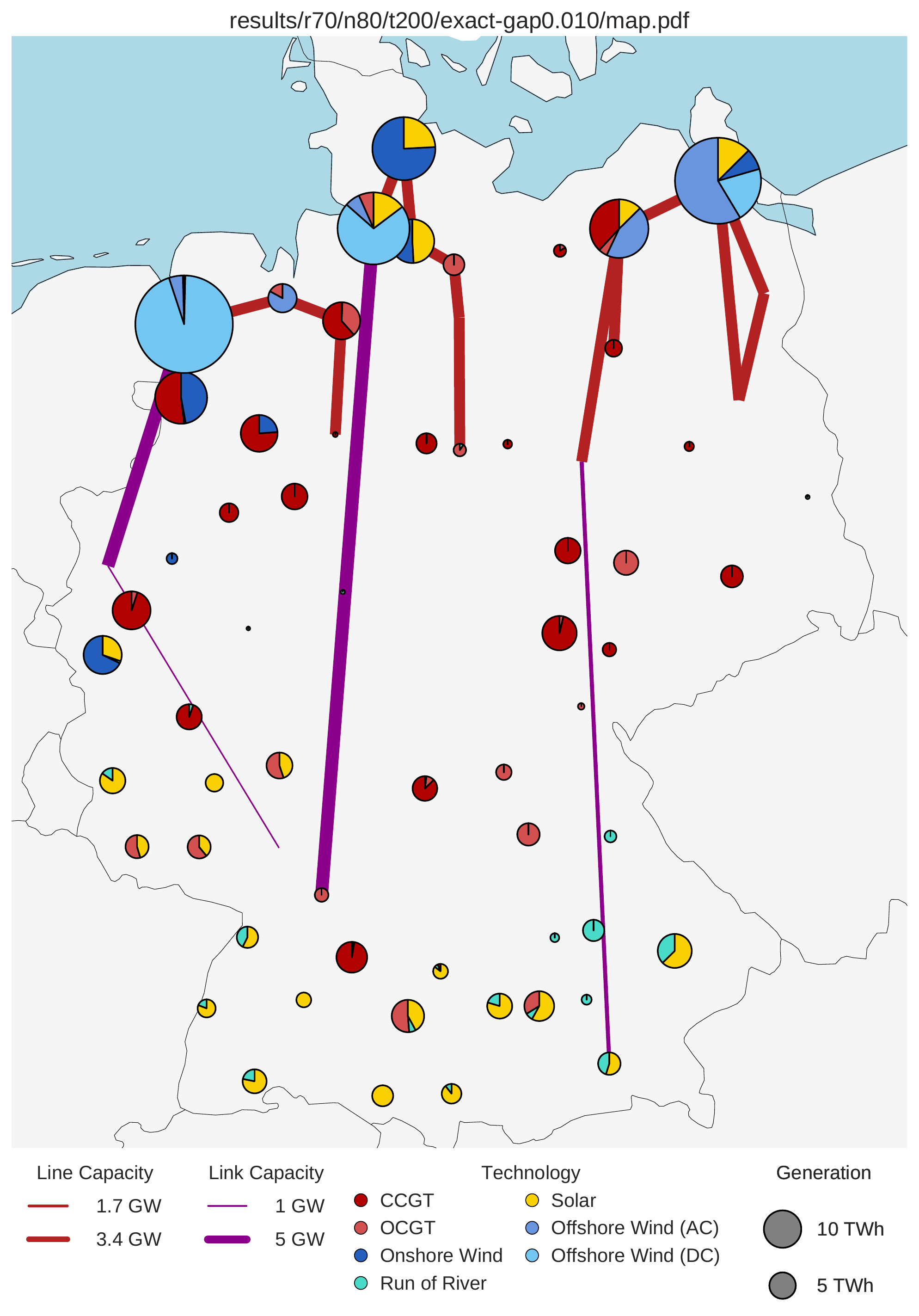} &
        \includegraphics[width=0.34\textwidth, trim=0cm 3.7cm 0cm 1cm, clip]{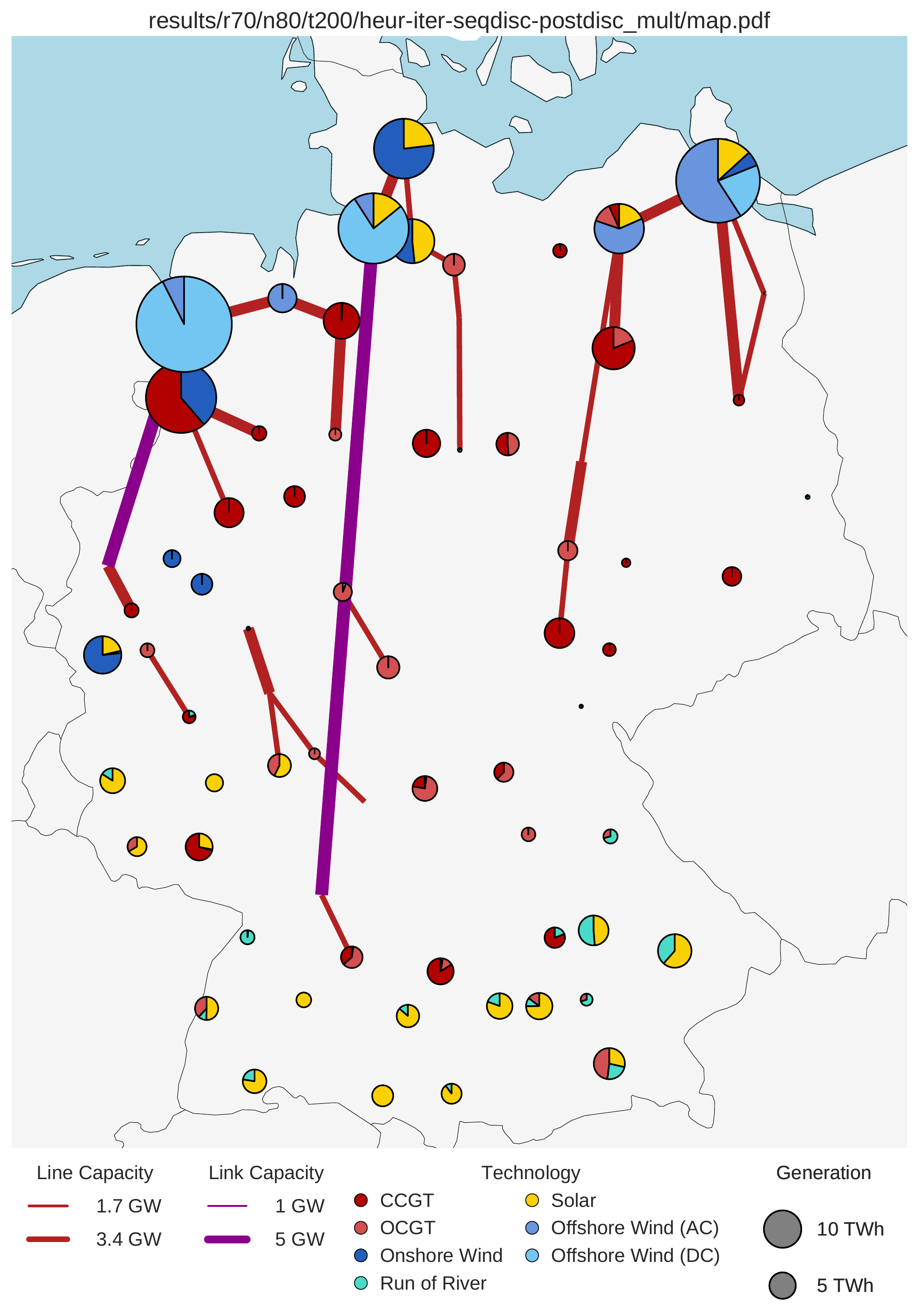} &
        \includegraphics[width=0.34\textwidth, trim=0cm 3.7cm 0cm 1cm, clip]{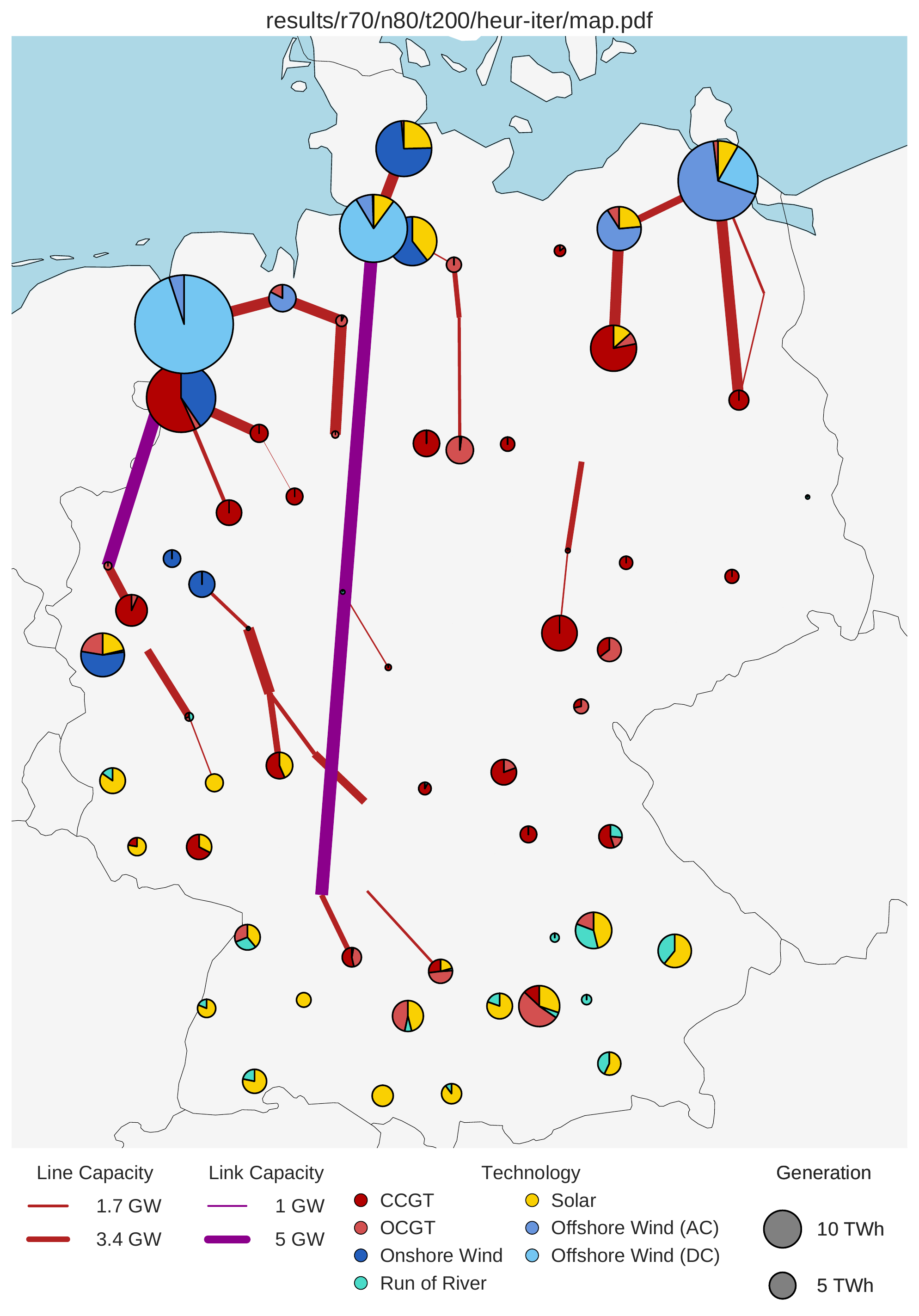} \\
      \end{tabular}
    \end{scriptsize}
  \end{center}
  \hfill \includegraphics[width=0.5\textwidth, trim=0cm 0cm 0cm 25.5cm, clip]{figures/map_n80_t200_cont.pdf}
  \vspace{-0.25cm} \caption{Map with transmission line expansion, HVDC link expansion and regional energy production for a model with 80 nodes and 200 snapshots obtained from
  (i) the exact big-$M$ formulation with a 1\% MIP optimality gap,
  (ii) an iterative heuristic approach with sequential impedance discretisation and post-facto line capacity discretization, and
  (iii) an iterative heuristic approach with continuous line expansion.
  }
  \label{fig:maps}
\end{figure*}

Figure \ref{fig:maps} offers a more visual approach to comparing the difference in transmission network expansion and power generation between the solutions of
(i) the big-$M$ formulation,
(ii) an iterative heuristic approach with sequential impedance discretisation and post-facto line capacity discretization, and
(iii) an iterative heuristic approach with continuous line expansion.

Across the three approaches, electricity production from renewables differs only marginally.
Likewise, although the capacities and production from OCGT and CCGT power plants may interchange between neighboring nodes, 
total volumes of gas-fired power plant capacities and production are very similar, as Table \ref{tab:gas} summarizes.
Thus, the regional generation mix is mostly insensitive to the applied heuristic.

\begin{table}
  \caption{Capacities and production of gas-fired powerplants for model and approaches of Figure \ref{fig:maps}}
  \begin{center}
    
    \begin{tabular}{llrrr}
      \hline
      & & \textbf{(i) big-$M$} & \textbf{(ii) integer} & \textbf{(iii) continuous} \\
      \hline
      GW  & CCGT & 21.6 & 22.2 & 21.4 \\
          & OCGT & 32.6 & 32.9 & 32.7 \\
      \hline
      TWh & CCGT & 96.9 & 96.1 & 96.8 \\
          & OCGT & 42.1 & 42.9 & 42.2 \\
      \hline
    \end{tabular}
  \end{center}
  \label{tab:gas}
\end{table}

Differences are more prevalent when looking at the expansion of HVDC links and AC transmission lines.
While the common theme is a focus of line expansion in Northern Germany, the heuristics prefer the expansion of individual lines in
Rhineland-Palatinate, Hesse and Saxony-Anhalt over building 2.5 GW on the route of the SuedOstLink (connecting the North-East of Germany with the South-East),
 which is the preferred choice of the big-$M$ formulation.
 Both, big-$M$ formulation and heuristics agree in utilizing 8 GW HVDC links on the routes of SuedLink and A-North
 (connecting the North-South axis in the West of Germany).


\begin{figure}
  \includegraphics[width=0.5\textwidth]{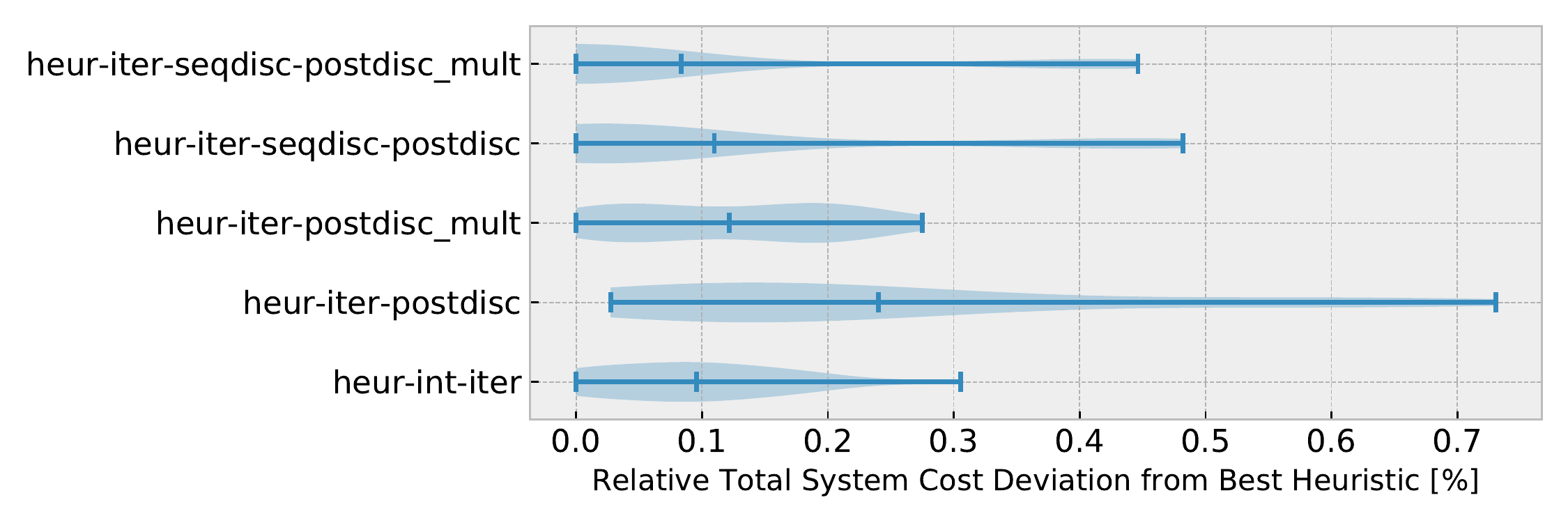}
  \vspace{-0.75cm}
  \caption{Comparison among heuristics in terms of relative total system cost.
  For each model, the objective values are related to the one heuristic approach with the lowest objective value.
  The distribution of relative total system cost deviations from the best heuristic is depicted here.}
  \label{fig:violincost}
\end{figure}

\begin{figure}
  \includegraphics[width=0.5\textwidth]{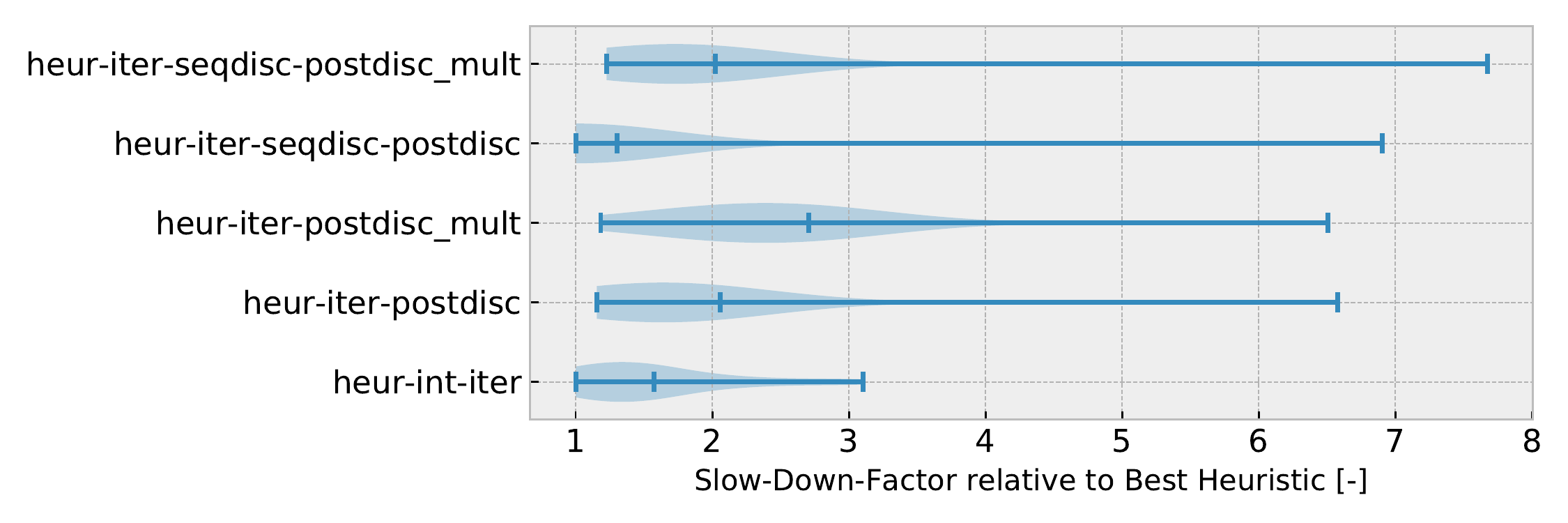}
  \vspace{-0.75cm}
  \caption{Comparison among heuristics in terms of solving times.
  For each model, the solving times are related to the fastest heuristic approach.
  The resulting distribution of the solving time deviations is depicted here.}
  \label{fig:violintime}
\end{figure}

Since the big-$M$ formulation was only solved for up to 200 snapshots, Figure \ref{fig:violincost} and Figure \ref{fig:violintime} extend the comparison between heuristics to up to 400 snapshots and address the quality of approaches as well as their solving times. On average all heuristics perform similarly well and, neglecting rare outliers, differ only within a range of 0.3\% in objective value.
A glance at the solving times reiterates that sequential discretization is an effective way of reducing the required number of iterations without a loss in solution quality. In contrast, evaluating multiple discretization thresholds consumes additional time which does not benefit solution quality by much.
Although the integer heuristic \textit{heur-int-iter} produces good quality results in comparably little time, it is unlikely to scale well for much larger problems due to its integer, non-convex nature.
Moreover, it is very sensitive to the choice of the  MIP optimality gap tolerance, which in this case was set to 0.5\%.
Additionally, solving times of integer problems can be very volatile. Such variance stems from randomness in the heuristics of most solvers (e.g. Gurobi).
Instead, the heuristic \textit{heur-iter-seqdisc-postdisc} seems to offer a sensible and well-scaling tradeoff between solving time and finding low-cost solutions. 


\begin{figure}
  \includegraphics[width=0.5\textwidth]{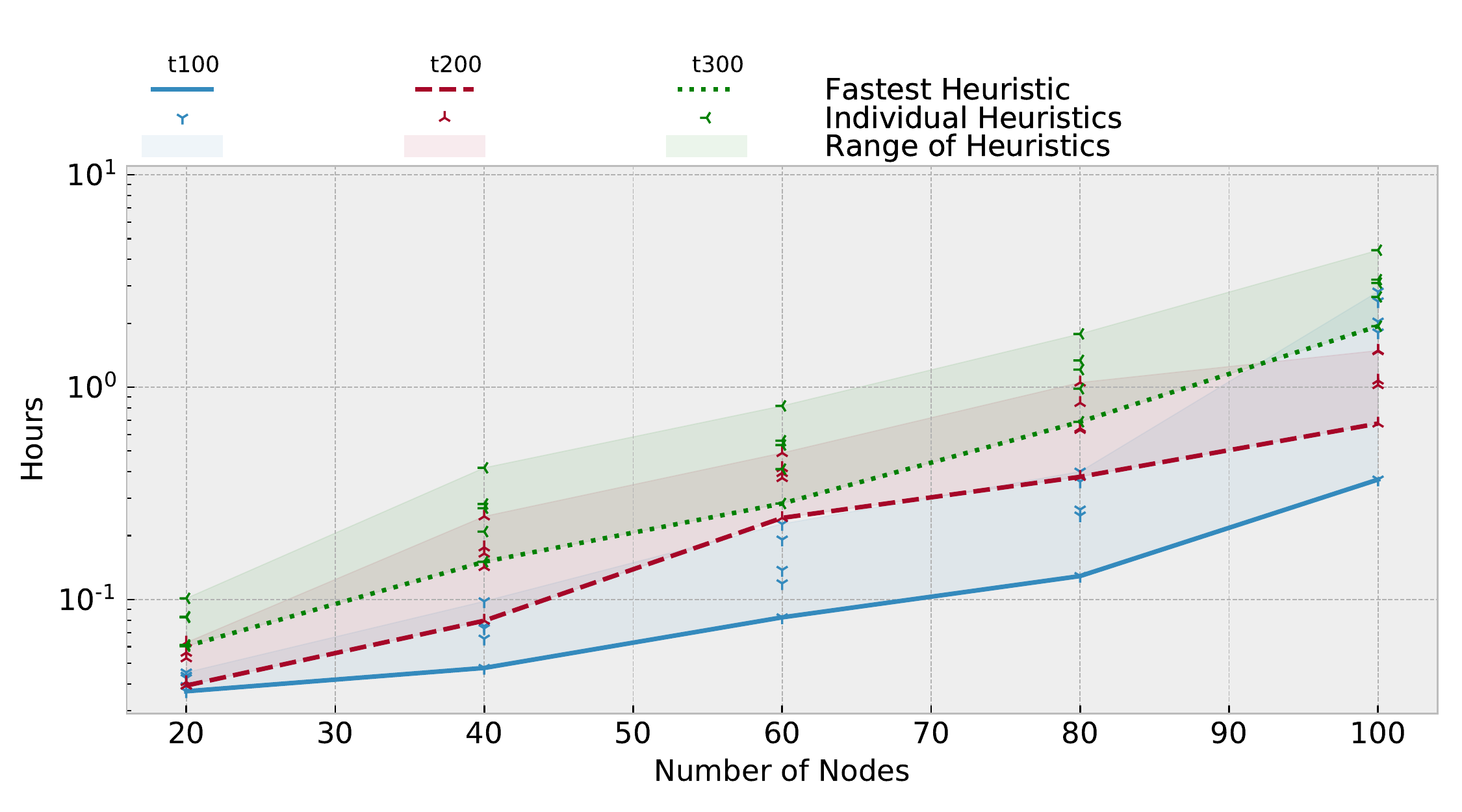}
  \vspace{-0.75cm} \caption{Scaling of solving time of heuristic approaches with the number of nodes for models with 100, 200, and 300 snapshots on a logarithmic scale.}
  \label{fig:scalingnodes}
\end{figure}

\begin{figure}
  \includegraphics[width=0.5\textwidth]{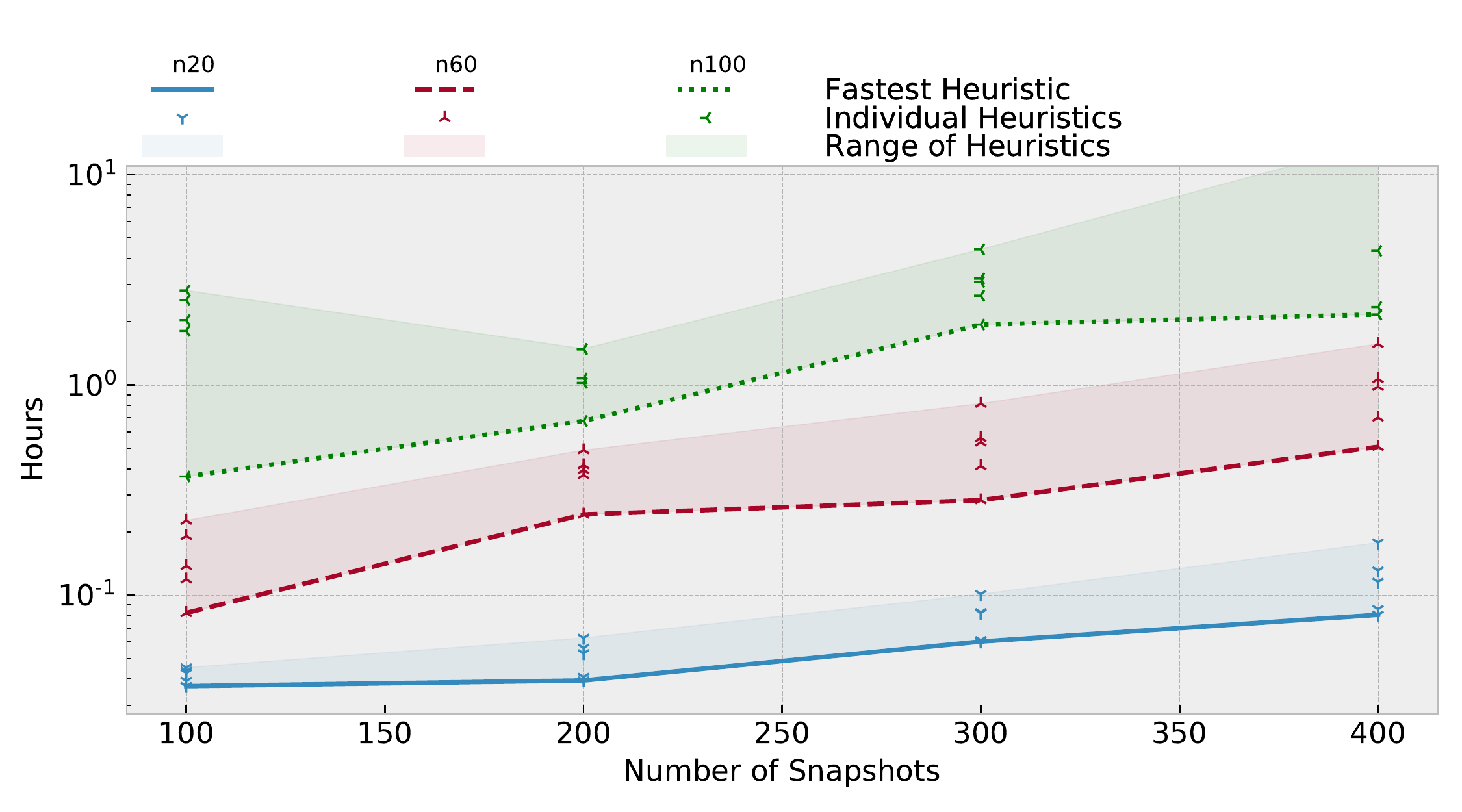}
  \vspace{-0.75cm} \caption{Scaling of solving time of heuristic approaches with the number of snapshots for models with 20, 60, 100 nodes on a logarithmic scale.}
  \label{fig:scalingsnapshots}
\end{figure}

To further investigate scaling, 
Figure \ref{fig:scalingnodes}
and
Figure \ref{fig:scalingsnapshots}
depict how much longer the heuristic approaches take when the number of nodes or the number of snapshots is increased from two different perspectives.
Both graphs show, that the sensitivity of solving times to model size is large and exponential, and that the heuristics seem to deal better with an increasing number of operating scenarios than with an increasing network size.


\begin{figure}
  \includegraphics[width=0.5\textwidth]{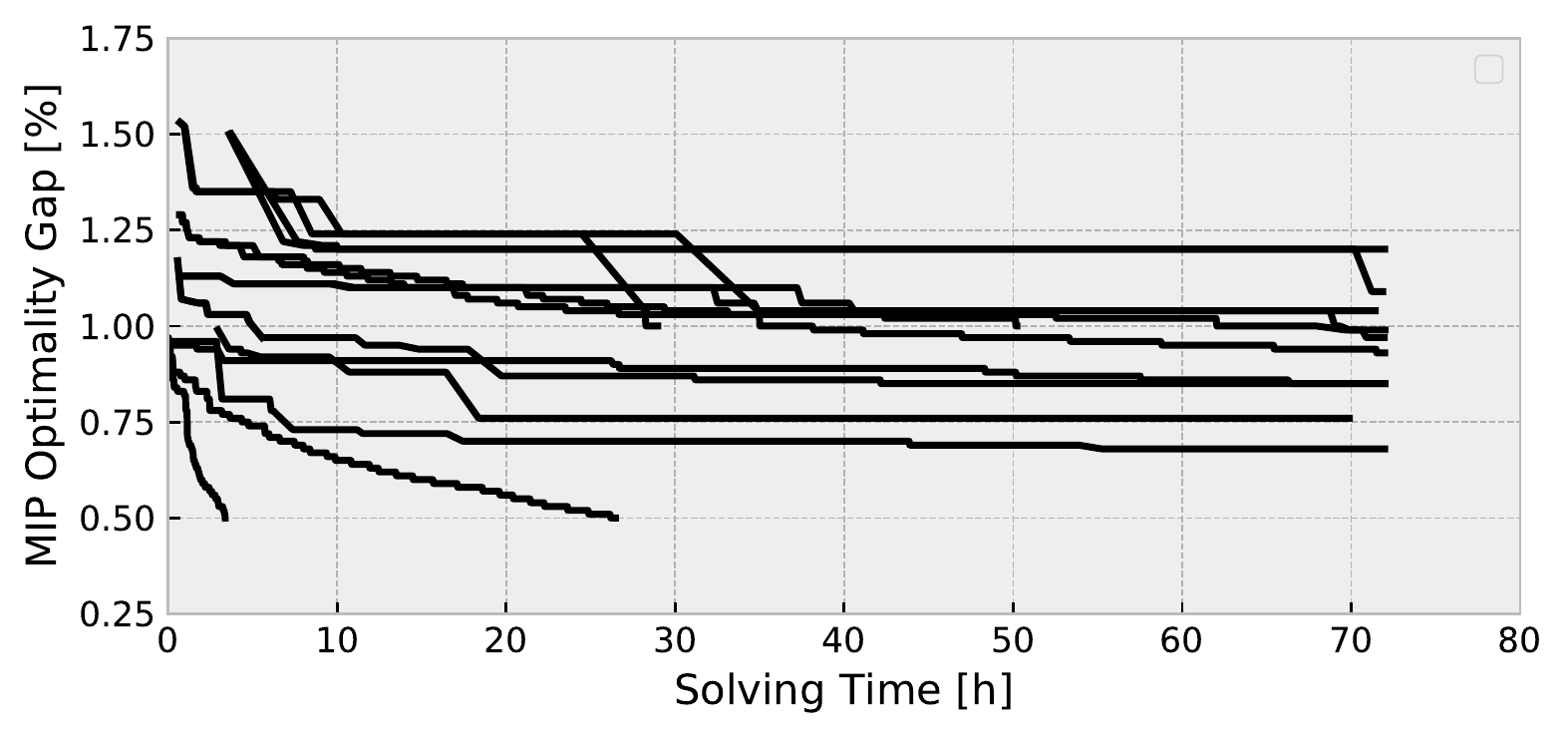}
  \vspace{-0.75cm} \caption{Progression of MIP optimality gaps among all model optimization runs which use the big-$M$ reformulation. Lines start when the first feasible solution is found and end when either the MIP optimality gap tolerance is satisfied or the walltime of the solver is reached.}
  \label{fig:gaps}
\end{figure}

As previously noted, when solving the big-$M$ formulation MILP, the upper bound is the best currently found feasible integer solution. 
Figure \ref{fig:gaps} shows, that a feasible solution which lies within a 1.5\% tolerance band can be found relatively quickly, but the progression from there becomes rapidly slower and may ultimately be terminated by the specified walltime.
In consequence, the benefit of large runtimes is questionable. However, it is unclear whether feasible solutions with similarly low initial optimality gaps can also be found for a full European model and hourly time resolution for a full year. 
We could however observe, that the extent of allowed line expansion significantly influences solving time for the big-$M$ formulation, and that the more ambitious the renewable target is chosen, the longer it takes to obtain an optimal solution.



\section{Critical Appraisal}
\label{sec:appraisal}

In this paper, the proximity of the objective functions for the heuristics are evaluated on relatively small networks and few snapshots,
which in general would not provide a sufficient basis for reliable investment decisions. But, it allows a reasonable performance comparison of the heuristic approaches to exact formulations.
Moreover, we disregarded storage options to avoid strong temporal-coupling and furthermore ignored the currently existing fleet of fossil-fueled power plants, but limited those to gas-fired options.
Neither do we consider sector-coupling, nor the transmission losses, unit commitment, or the provision of ancillary services.
The test models also exclusively cover Germany, while questions of transmission expansion are in reality closely intertwined with considerations of other European countries.

Besides remedying these issues, future work further involves analysing the solution space.
The fact that the different heuristics provide solutions very close to the global optimum but differ in their individual investment decisions,
indicates that many near-optimal solutions exist which may be less susceptible to public or regulatory opposition.

\section{Conclusions}
\label{sec:conclusion}
 
In this paper we have compared several heuristics for approximating the exact solution
of joint transmission and generation expansion planning
with the linear approximation of the power flow equations.

Particularly in light of the complexity already a limited choice of reinforcements entails, we conclude
for models with high temporal and spatial resolution 
that using a continuous relaxation of line investment together with subsequent post-facto discretization and
applying sequential discretization of line impedances to accelerate convergence is beneficial.
This particular heuristic could reduce computation times of the joint optimisation by 82\% with a maximal total system cost deviation of 1.5\% for the models considered here.

Altogether, it was shown that the presented heuristics closely mirror optimal integer line investment of the exact MINLP with 
considerable savings in solving time for policy-relevant low-carbon energy system optimization models.
With larger models that
include storage and sector coupling it can be expected that only the
heuristic methods will solve within a practical finite time.

\section*{Acknowledgments}

F.N. and T.B. gratefully acknowledge funding from the Helmholtz
Association under grant no. VH-NG-1352.



\bibliographystyle{IEEEtran}
\bibliography{library.bib}
%



\end{document}